\input amstex
\documentstyle{amsppt}
\input bull-ppt
\keyedby{bull503e/paz}


\def\crossymbol{\hbox{%
    \hbox{$\,$\vrule height 10 pt depth 7pt width 2.5 pt}\kern-6.6pt%
    \hbox{\raise4pt\hbox{\vrule height 2.5 pt depth 0pt width 11pt}$\,$}}}

\long\def\block #1\endblock{\vskip 6pt
	{\leftskip=3pc \rightskip=3pc
	\noindent #1\endgraf}\vskip 6pt}

\catcode`\@=11
\def\raggedcenter{\leftskip\z@ plus.4\hsize 
\rightskip\leftskip
 \parfillskip\z@ \parindent\z@ \spaceskip.3333em 
\xspaceskip.5em
 \pretolerance9999\tolerance9999 \exhyphenpenalty\@M
 \hyphenpenalty\@M}
\catcode`\@=12

\def\center{\bgroup\raggedcenter}
\def\endcenter{\endgraf\egroup}

\def\lethead#1{\vskip6pt{\noindent\bf#1\par\nobreak%
\vskip6pt}}

\def\paragraph#1{\vskip6pt\noindent{\bf#1}\par}

\def\letaddress{\bgroup\interlinepenalty=10000
	\bf \parindent0pt\parskip0pt\obeylines}
\def\endletaddress{\egroup\nobreak\vskip1pc}

\parindent18pt
\parskip0pt
\overfullrule0pt

\def\newletter{\vskip3pc}

\let\sl\it


\topmatter
\cvol{30}
\cvolyear{1994}
\cmonth{April}
\cyear{1994}
\cvolno{2}
\cpgs{178-207}
\title Responses to ``Theoretical Mathematics: Toward\\
a cultural synthesis of mathematics and\\
theoretical physics'', by A. Jaffe and F. Quinn\endtitle
\shorttitle{Responses to ``Theoretical Mathematics...''}
\author Michael Atiyah et al.\endauthor
\subjclass Primary 01A80\endsubjclass
\endtopmatter

\document

\letaddress
Michael Atiyah
The Master's Lodge
Trinity College
Cambridge CB2 1TQ
England, U.K.
\endletaddress

 	I find myself agreeing with much of the detail of the 
Jaffe-Quinn
argument, especially the importance of distinguishing 
between results based on
rigorous proofs and those which have a heuristic basis.  
Overall, however, I
rebel against their general tone and attitude which 
appears too authoritarian.

	My fundamental objection is that Jaffe and Quinn present 
a sanitized
view of mathematics which condemns the subject to an 
arthritic old age.  They
see an inexorable increase in standards of rigour and are 
embarrassed by
earlier periods of sloppy reasoning.  But if mathematics 
is to rejuvenate
itself and break exciting new ground it will have to allow 
for the exploration
of new ideas and techniques which, in their creative 
phase, are likely to be as
dubious as in some of the great eras of the past.  Perhaps 
we now have high
standards of proof to aim at but, in the early stages of 
new developments, we
must be prepared to act in more buccaneering style.

	The history of mathematics is full of instances of happy 
inspiration
triumphing over a lack of rigour.  Euler's use of wildly 
divergent series or
Ramanujan's insights are among the more obvious, and 
mathematics would have
been poorer if the Jaffe-Quinn view had prevailed at the 
time.  The marvelous
formulae emerging at present from heuristic physical 
arguments are the modern
counterparts of Euler and Ramanujan, and they should be 
accepted in the same
spirit of gratitude tempered with caution.

	In fact the whole area between Quantum Field Theory and 
Geometry (which
is the main target of Jaffe-Quinn) has now produced a 
wealth of new results
which have striking evidence in their favour. In many 
important cases we now
have rigorous proofs based on other methods. This provides 
additional
confidence in the heuristic arguments used to discover the 
results in the first
place.

	For example, Witten's work has greatly extended the scope 
of the Jones
knot invariants and, when the dust settles, I think we 
will see here a fully
rigorous topological quantum field theory in 2+l 
dimensions.  The Feynman
integrals will have been given precise meanings, not by 
analysis, but by a
mixture of combinatorial and algebraic techniques.  The 
disparaging remarks in
the Jaffe-Quinn article are totally unjustified.

	Other results of this type (which now have rigorous 
proofs) include the
rigidity of the elliptic genus, formulae for the volume of 
moduli spaces,
computations of their cohomology, and information about 
rational curves on 
3-dimensional Calabi-Yau manifolds.  A new and much 
simpler proof of the positive
energy theorem of Schoen and Yau emerged from ideas of 
Witten, based on the
Dirac operator and a super-symmetric formalism.  
Jaffe-Quinn single out the
Schoen-Yau proof as an example of respectable mathematical 
physics, but they
deny the title to Witten.

	While acknowledging the important role of conjectures in 
mathematics,
Jaffe and Quinn reserve their garland for the person who 
ultimately produces
the rigorous proof.  For example, they cite the famous 
Weil conjectures and the
eventual proof by Deligne (and Grothendieck).  But surely 
Weil deserves
considerable credit for the whole conception (and the 
proof for curves)?  The
credit which posterity ascribes depends on the respective 
weight of ideas and
techniques in the conjecture and its proof.  In the case 
of Hodge's theory of
harmonic forms, Hodge's own proof was essentially faulty 
because his
understanding of the necessary analysis was inadequate. 
Correct proofs were
subsequently provided by better analysts, but this did not 
detract from Hodge's
glory. The mathematical world judged that Hodge's 
conceptual insight more than
compensated for a technical inadequacy.

	Jaffe represents the school of mathematical physicists 
who view their
role as providing rigorous proofs for the doubtful 
practices of physicists.
This is a commendable objective with a distinquished 
history.  However, it
rarely excites physicists who are exploring the front line 
of their subject.
What mathematicians can rigorously prove is rarely a hot 
topic in physics.

	What is unusual about the current interaction is that it 
involves
front-line ideas both in theoretical physics and in 
geometry.  This greatly
increases its interest to both parties, but Jaffe-Quinn 
want to emphasize the
dangers.  They point out that geometers are inexperienced 
in dealing with
physicists and are perhaps being led astray. I think most 
geometers find this
attitude a little patronizing: we feel we are perfectly 
capable of defending
our virtue.

	What we are now witnessing on the geometry/physics 
frontier is, in my
opinion, one of the most refreshing events in the 
mathematics of the 20th
century.  The ramifications are vast and the ultimate 
nature and scope of what
is being developed can barely be glimpsed.  It might well 
come to dominate the
mathematics of the 2lst century.  No wonder the younger 
generation is being
attracted, but Jaffe and Quinn are right to issue warning 
signs to potential
students.  For those who are looking for a solid, safe PhD 
thesis, this field
is hazardous, but for those who want excitement and action 
it must be
irresistible.

\newletter


\letaddress
Armand Borel
Institute for Advanced Study
Princeton, NJ 08540
borel@math.ias.edu
\endletaddress

Some comments on the article by A. Jaffe and F. Quinn:

There are a number of points with which I agree, but they 
are so
obvious that they do not seem worth such an elaborate 
discussion.
On the other hand, I disagree with the initial stand and 
with much
of the general thrust of the paper, so that I shall not 
comment
on it item by item, but limit myself to some general 
remarks.

First the starting point.  I have often maintained, and even
committed to paper on some occasions, the view that 
mathematics is
a science, which, in analogy with physics, has an 
experimental
and a theoretical side, but operates in an intellectual 
world of
objects, concepts and tools.
Roughly, the experimental side is the investigation of 
special
cases, either because they are of interest in themselves or
because one hopes to get a clue to general phenomena, and 
the
theoretical side is the search of general theorems.  In 
both, I
expect proofs of course, and I reject categorically a 
division
into two parts, one with proofs, the other without.

I also feel that what mathematics needs least are pundits 
who
issue prescriptions or guidelines for presumably
less enlightened mortals.
Warnings about the dangers of certain directions are of 
course
nothing new.
In the late forties, H. Weyl was very worried by the trend
towards abstraction, exemplified by the books of Bourbaki 
or that
of Chevalley on Lie groups, as I knew from M. Plancherel.
Later, another mathematician told me he had heard such views
from H. Weyl in the late forties but, then, around 1952 I
believe, i.e. after the so-called French explosion, H. Weyl
told him:  ``I take it
all back.''

In fact, during the next quarter century, we experienced a
tremendous development of pure mathematics, bringing 
solutions
of one fundamental problem after the other, unifications, 
etc.,
but during all that time, there was in some quarters some 
whining
about the dangers of the
separation between pure math and applications to sciences, 
and
how the great nineteenth century mathematicians cultivated 
both
(conveniently ignoring some statements by none other than 
Gauss
which hardly support that philosophy).
[To avoid any misunderstanding, let me hasten to add that 
I am
not advocating the separation between the two, being quite
aware of the great benefits on both sides of interaction, 
but
only the freedom to devote oneself to pure mathematics, if 
so
inclined.]

Of course, I agree that no part of mathematics can 
flourish in
a lasting way without solid foundations and proofs, and 
that not
doing so was harmful to Italian algebraic geometry for 
instance.
I also feel that it is probably so for the Thurston 
program, too.
It can also happen that standards of rigor deemed 
acceptable by
the practitioners in a certain area turn out to be found 
wanting
by a greater mathematical community.  A case in point, in my
experience, was E. Cartan's work on exterior differential 
forms
and connections, some of which was the source of a rather 
sharp
exchange between Cartan and Weyl.
Personally, I felt rather comfortable with it but later, 
after
having been exposed to the present points of view, could 
hardly 
understand what I had thought to understand.  We all know 
about 
Dirac diagonalizing any self-adjoint operator and using 
the Dirac
function.  And there are of course many more examples.  
But I do
believe in the self-correcting power of mathematics, already
expressed by D. Hilbert in his 1900 address, and all I have
mentioned (except for Thurston's program) has been 
straightened
out in due course.  Let me give another example of this
self-correcting power of mathematics.  In the early 
fifties, the
French explosion in topology  was really algebraic 
topology with
a vengeance.
Around 1956, I felt that topology as a whole
was going too far in that direction and I
was wishing that some people would again get their hands 
dirty
by using more intuitive or geometric points of view (I did
so for
instance in a conversation with J. C. Whitehead at the time,
 which he
reminded me of shortly before his death, in 1960; I had
forgotten it.)
But shortly after came the developments of $PL$-topology
by Zeeman and Stallings, of differential topology by Milnor
and Smale, and there was subsequently in topology a 
beautiful
equilibrium between algebraic, differential and $PL$ 
points of
view.

But this was achieved just because some gifted people 
followed their own
inclinations, not because they were taking heed of some 
solemn
warning.

In advocating freedom for  mathematicians, I am not 
innovating
at all.  I can for instance refer to a lecture by A.~Weil
({\it Collected Papers\/} II, 465--469) praising 
disorganization
in mathematics and pointing out that was very much the
way Bourbaki operated.
As a former member of Bourbaki, I was of course saddened to
read that all that collective work, organized or not,
ended up with the erection of a bastion
of arch-conservatism.
Not entertaining pyramids of conjectures?
Let me add that Weil was not ostracized for his conjectures,
nor was Grothendieck for his standard conjectures and the
theory motives, nor Serre for his ``questions''.

F. Quinn is not making history in raising questions about 
the
Research Announcements in the {\it Bulletin\/}, as you 
know.  At some
point, they were functioning poorly and their
suppression was suggested by some.
To which I. Singer answered that people making such 
proposals
did not know what the AMS was about (or something to that
effect) and offered to manage that department for a few 
years.
He did so and it functioned very well during his
tenure.
Also, the {\it Comptes Rendus\/} have a very long history 
of R.A.s.
There were
ups and downs of course, but for the last twenty years
or so, it seems to me to have been working well on the 
whole.
All this to say that the problems seen by F. Quinn are not
new, have been
essentially taken care of in the past,
and I do not see the need for new
prescriptions.

\newletter


\letaddress
G. J. Chaitin
IBM Research Division
P. O. Box 704
Yorktown Heights, NY 10598
chaitin@watson.ibm.com
\endletaddress

\block
{\bf Abstract.} It is argued that the information-theoretic
incompleteness theorems of algorithmic information theory 
provide a
certain amount of support for what Jaffe and Quinn call 
``theoretical
mathematics''.
\endblock

One normally thinks that everything that is true is true 
for a reason.
I've found mathematical truths that are true for no reason 
at all.
These mathematical truths are beyond the power of 
mathematical
reasoning because they are accidental and random.

Using software written in {\sl Mathematica\/} that runs on 
an IBM
RS/6000 workstation [5, 7], I constructed a perverse 
200-page exponential
diophantine equation with a parameter $N$ and 17,000 
unknowns:

\center
        Left-Hand-Side($N$) = Right-Hand-Side($N$).
\endcenter

For each nonnegative value of the parameter $N$, ask 
whether this
equation has a finite or an infinite number of nonnegative 
solutions.
The answers escape the power of mathematical reason 
because they are
completely random and accidental.

This work is part of a new field that I call algorithmic 
information
theory [2,3,4].

What does this have to do with Jaffe and Quinn [1]?

The result presented above is an example of how my
information-theoretic approach to incompleteness makes 
incompleteness
appear pervasive and natural.  This is because algorithmic 
information
theory sometimes enables one to measure the information 
content of a
set of axioms and of a theorem and to deduce that the 
theorem cannot
be obtained from the axioms because it contains too much 
information.

This suggests to me that sometimes to prove more one must 
assume more,
in other words, that sometimes one must put more in to get 
more out.
I therefore believe that elementary number theory should 
be pursued
somewhat more in the spirit of experimental science.  
Euclid declared
that an axiom is a self-evident truth, but physicists are 
willing to
assume new principles like the Schr\"odinger equation that 
are not
self-evident because they are extremely useful.  Perhaps 
number
theorists, even when they are doing elementary number 
theory, should
behave a little more like physicists do and should 
sometimes adopt new
axioms.  I have argued this at greater length in a lecture 
[6, 8] that I
gave at Cambridge University and at the University of New 
Mexico.

In summary, I believe that the information-theoretic 
incompleteness
theorems of algorithmic information theory [2,3,4,5,6,7,8] 
provide a
certain amount of support for what Jaffe and Quinn [1] call
``theoretical mathematics''.

\lethead{References}

\item {[1]}  A. Jaffe and F. Quinn, {\it Theoretical 
Mathematics\,\RM: Toward
a cultural synthesis of mathematics and theoretical 
physics\/},
Bull. Amer. Math. Soc. {\bf 29} (1993),
1-13.
\item {[2]} G. J. Chaitin, {\it Algorithmic information 
theory,} revised third printing,
Cambridge Univ. Press, Cambridge, 1990.
\item {[3]} G. J. Chaitin, {\it Information, randomness \& 
incompleteness---%
papers on algorithmic information theory,} second edition,
World Scientific, Singapore, 1990.
\item {[4]} G. J. Chaitin, {\it Information-theoretic 
incompleteness,}
World Scientific, Singapore, 1992.
\item {[5]} G. J. Chaitin,
{\it Exhibiting randomness in arithmetic using Mathematica 
and C\/},
IBM Research Report RC-18946, June 1993.
\item {[6]} G. J. Chaitin,
{\it Randomness in arithmetic and the decline and fall of 
reductionism in
pure mathematics\/}, {Bull. European Assoc. for Theoret. 
Comput. Sci.,} no.\ 50 (July 1993). 
\item {[7]} G. J. Chaitin,
{\it The limits of mathematics---Course outline and 
software\/}, IBM Research
Report RC-19324, December 1993.
\item {[8]} G. J. Chaitin, {\it Randomness and complexity 
in pure
mathematics\/}, Internat. J. Bifur. Chaos {\bf4} (1994) 
(to appear).

\vfill


\newletter

\letaddress
Daniel Friedan 
Department of Physics
Rutgers University
New Brunswick, NJ  08903
friedan@physics.rutgers.edu
\endletaddress

The paper distorts the relation of experiment to 
theoretical physics.
To paraphrase Fermi (perhaps badly): an experiment which 
finds the
unexpected is a discovery; an experiment which finds the 
expected is a
measurement.

I have the impression that applying rigor to a theoretical 
idea is
given substantial credit when it disconfirms the 
theoretical idea or
when the proof is especially difficult or when the ideas 
of the proof
are original, interesting and fruitful.  This seems quite 
enough to
motivate the application of rigor, for those who are 
motivated by the
prospect of credit.  Perhaps pedestrian proofs do get only 
a little
recognition, but should they really get more?
Is it useful to formulate explicit general rules for 
assigning credit
in mathematics?

Is there really any evidence that mathematics is suffering 
from the
theoretical influence?  Are mathematicians really finding 
it difficult
to read theoretical papers critically, detecting for 
themselves the
level of rigor?  Are rigorous-minded graduate students so 
awash in
problems that they truly resent the offerings of the 
so-called
theoretical mathematicians?

As far as I know, there has never been a surplus of 
originality in
mathematics or in physics.  Is it useful to criticize the 
manner of
expression of original ideas on the grounds that the 
community is slow
to absorb, evaluate and/or pursue them?


\newletter

\letaddress
James Glimm
Department of Applied Mathematics and Statistics
State University of New York at Stony Brook
Stony Brook, NY 11794-3600
mills@ams.sunysb.edu
\endletaddress

Truth, in science, lies not in the eye of the beholder, 
but in objective
reality.  It is thus reproducible across barriers of 
distance, political
boundaries and time.  As mathematics becomes increasingly 
involved in
interdisciplinary activities, clashes with distinct 
standards of proof
from other disciplines will
arise.  The Jaffe-Quinn article is thus constructive, in 
opening and framing
this discussion for the interaction of mathematics with 
physics.

These issues are older, and perhaps better understood, 
within the
applied mathematics community.  The outcome there follows 
the broad outlines
proposed by Jaffe and Quinn: clear labeling of standards 
which are adopted
within a specific paper, or ``truth in advertising''.   
Additionally,
especially for computational mathematics, the standard of 
reproducibility
is tested by the actual reproduction of a (similar, 
related or even
identical) experiment, by (say) other methods.  However, 
the most central
standard of truth in science is the agreement between 
theory and data
(e.g. laboratory experiments).  In this sense, science has 
a standard which
goes beyond that of mathematics.  A conclusion is correct 
according to
the standards of science if \underbar{both} the hypotheses 
and the
reasoning connecting the hypotheses to the conclusion are 
valid.

It bears repeating that the correct standards for 
interdisciplinary work
consist not of the intersection, but the union of the 
standards from the two
disciplines. Specificially, speculative theoretical 
reasoning in physics is
usually strongly constrained by experimental data.  If 
mathematics is going to
contemplate a serious expansion in the amount of 
speculation which it supports
(which could have positive consequences), it will have a 
serious and
complementary need for the admission of new objective 
sources of data, going
beyond rigorously proven theorems, and including computer 
experiments,
laboratory experiments and field data.  Put differently, 
the absolute standard
of logically correct reasoning was developed and tested in 
the crucible of
history.  This standard is a unique contribution of 
mathematics to the culture
of science.  We should be careful to preserve it, even (or 
especially) while
expanding our horizons.


\newletter

\letaddress
Jeremy J. Gray 
Faculty of Mathematics
Open University
Milton Keynes, MK7 6AA
England
\endletaddress

	  The letter by Jaffe and Quinn raises several issues.  
The extent to
which theoretical mathematics, as they term it, is 
prevalent in mathematics is
perhaps for a mathematician rather than an historian of 
mathematics to comment
upon, but it seems to me that one aspect of the problem is 
underestimated.  Not
only students and young researchers but all those who work 
away from the main
centers of research are disadvantaged.  They too will be 
encouraged to rely
unwisely on insecure claims.  They will also be unaware of 
the degree to which
the claims of theoretical mathematics are discounted or 
interpreted by experts
in the field.

	Indeed, the role of experts in this connection is more 
complicated than
the authors have suggested. Preprints aside, theoretical 
mathematics is
published presumably because competent referees have 
endorsed it. It is often
accompanied by talks, invited lectures, and conference 
papers given not only by
the author but others equally convinced of the merits of 
the work. The problem
does not arise merely with one mathematician claiming too 
much, but with a
network of others endorsing the claims.

	The involvement of experts points to a problem with the 
remedies
proposed by Jaffe and Quinn.  Their plea for honest 
advertising is surely to be
accepted.  The difficult problem we all have to confront 
is the honest mistake.
What of the paper that offers a result which is not, in 
the opinion of the
author, conjectural but proved?  To be sure, the proof is 
based on insights
that have not yet yielded to expression in the form of 
definitions, lemmas, and
proofs, but it seems clear to the author. It may well be 
clear to the referee.
Both would wish to see the work presented as rigorous 
mathematics. Many would
argue that the tricky concept of shared insight rather 
than logical precision
is what mathematical communication is about. But since 
there are no absolute
canons of rigour, and it is impossible to insist that 
every paper be written so
that a (remarkably) patient graduate student can follow 
it, some mistakes are
inevitably published.  This observation does not render 
the proposed remedy
nugatory, but it suggests that we shall still be working 
in an imperfect world.

	Care would also have to be given to the suggestion that 
theoretical
work could be published as such.  The present system has 
the virtue that
discoverers of new and important mathematics work as hard 
as they can to prove
the validity of their claims.  So far as I know, the 
cautionary tales presented
by Jaffe and Quinn are all tales of mathematicians who at 
the time of
publication had done their best to present correct 
statements.  On the other
hand, it is well known that working mathematicians 
habitually entertain ideas
which do not quite work out as they had hoped.  Their 
initial insight was not,
after all, veridical.  One would have to be cautious of 
adopting a scheme
whereby a good mathematician could publish a paper, 
labelled theoretical,
without trying flat out to prove its results. There are, 
at the highest level,
few if any more likely to come up with the proofs than 
these creative
mathematicians themselves.

	The historical examples given are also open to 
refinement, but in ways
that if anything support the paper.  It is true that 
classic Italian algebraic
geometry entered a decline, and that by modern standards 
it seems to lack
rigour---but this perception is modern, and due to 
Zariski, who also brought
new questions to bear (such as arbitrary fields).  What it 
seemed to contain at
the time was a rich mixture of results and problems.  
Certain key topics were
held to be securely established at one moment, more 
doubtful at another, much
as is the case in some topics today.  Poincar\'e typically 
wrote papers that
few could respond to for a generation.  The reasons are 
not clear, but the
theoretical nature of his work cannot have helped.  A 
student of mine, June
Barrow-Green, has recently shown that a major mistake in 
his prize-winning
essay on celestial mechanics eluded the judges, 
Weierstrass and Hermite.  Had
Poincar\'e not spotted the error himself, it would 
presumably have been
published with their implicit endorsement.

	What is perhaps the greatest change over the last one 
hundred years is
not that standards have risen---the authors make a good 
case that they have
not---but that the profession has grown.  Poincar\'e may 
have had an audience
of no more than ten capable of following him at his most 
inspired, and they all
had consuming interests of their own.  Today's leading 
figures soon attract
seminars in half a dozen places, the attention of many 
other mature
mathematicians, and perhaps the zeal of ambitious graduate 
students.  What is
striking is that despite all this attention, there are 
still the problems to
which Jaffe and Quinn have drawn our attention and for 
which we must surely
thank them for outlining remedies.  Because the best 
theoretical work may
convince even experts unduly, I fear that we cannot be 
optimistic about the
outcome.

\vfill


\newletter

\letaddress
Morris W. Hirsch
Department of Mathematics
University of California at Berkeley
Berkeley, CA 94720-0001
hirsch@math.berkeley.edu
\endletaddress

\lethead{Theoretical, Speculative and Nonrigorous 
Mathematics}

Several interesting and controversial points are raised in 
this
provocative essay.

To begin with, the authors make up a new term, ``theoretical
mathematics''.  They suggest that there is a growing 
branch of mathematics called theoretical mathematics, 
whose relation
to rigorous mathematics is parallel to that between 
theoretical
physics and experimental physics.  They warn of dangers in 
this kind
of division of labor, but suggest that this new field 
could be a respectable
branch of mathematics.

 Even though the authors ``do not wish to get involved in 
a discussion
of terminology'', it is important to note at the outset 
that their use
of ``theoretical'' is tied to a controversial 
philosophical position:
that mathematics is about the ``nature of reality'',  
later qualified
as ``mathematical reality'', apparently distinct from 
``physical
reality''.  They suggest ``Mathematicians may have even 
better access
to mathematical reality than the laboratory sciences have 
to physical
reality.''  

While they wisely don't attempt to define ``mathematical
reality'', this philosophical stance complicates and 
prejudices the
discussion.  For if we don't assume that mathematical 
speculations are
about ``reality'' then the analogy with physics is greatly 
weakened---and there is then no reason to suggest that a 
speculative mathematical
argument is a theory of anything, any more than a poem or 
novel is
``theoretical''.  For this reason I hope this use of 
``theoretical''
is not generally adopted; instead I prefer the more 
natural {\it
speculative mathematics}.  It should be obvious that there 
is a huge
difference between theoretical physics and speculative 
mathematics!

The nonrigorous use of mathematics by scientists, 
engineers, applied
mathematicians and others, out of which rigorous 
mathematics sometimes
develops, is in fact more complex than simple speculation. 
 While
sloppy proofs are all too common, deliberate presentation 
of unproved
results as correct is fortunately rare.

Much more frequent is the use of mathematics for {\it 
narrative
purposes}.  An author with  a story to tell  feels it can
be expressed most clearly in mathematical language.  In 
order to tell
it coherently without the possibly infinite delay rigor 
might require,
the author introduces certain assumptions, speculations 
and leaps of
faith, e.g.: ``In order to proceed further we assume the 
series
converges---the random variables are independent---the 
equilibrium is
stable---the determinant is nonzero---.'' In such cases it 
is often
irrelevant whether the mathematics can be rigorized, 
because the
author's goal is to persuade the reader of the 
plausibility or
relevance of a certain view about how some real world 
system behaves.
The mathematics is a language filled with subtle and 
useful metaphors.
The validation is to come from experiment---very possibly 
on a
computer.  The goal in fact may be to suggest a particular 
experiment.
The result of the narrative will be not new mathematics, 
but a new
description of ``reality'' ({\it real} reality!) 

This use of mathematics can be shocking to the pure 
mathematician
encountering it for the first time; but it is not only 
harmless, but
indispensable to scientists and engineers.

\lethead{Was Poincar\'{e} Speculative?}
We must carefully distinguish between modern papers 
containing
mathematical speculations, and papers published a hundred 
years ago
which we, today, consider defective in rigor, but which 
were perfectly
rigorous according to the standards of the time.  
Poincar\'{e} in his
work on Analysis Situs was being as rigorous as he could, 
and
certainly was not consciously speculative.  I have seen no 
evidence
that contemporary mathematicicans considered it 
``reckless'' or
``excessively theoretical''.  When young Heegard in his 1898
dissertation brashly called the master's attention to 
subtle mistakes,
Poincar\'{e} in 1899, calling Heegard's paper ``tr\`es 
remarquable'',
respectfully admitted his errors and repaired them.  In 
contrast, in
his 1912 paper on the Annulus Twist theorem (later proved by
Birkhoff), Poincar\'{e} apologized for publishing a 
conjecture, citing
age as his excuse.

I don't accept the authors' Cautionary Tale of the ``slow 
start'' in
algebraic and differential topology due to Poincar\'{e}'s 
having
``claimed too much, and proved too little''.  In fact he 
proved quite
a lot by the standards of the day---but there was little 
use for it
because the mathematics which could use it was not 
sufficiently
developed.  The ``15 or 20 years'' which it took for ``real
development to begin'' was not a long period in that more 
leisurely
age.  
In fact it was not until the late 1950s that what is now 
called
differential topology found substantial application in 
other fields.

Poincar\'e could indeed be careless: In 1900 he announced 
that if (in
our terminology) a closed manifold has the same Betti 
numbers as the
3-sphere $S^3$ then it is homeomorphic to $S^3$.  But in 
1904 he
admitted his error (neglecting the fundamental group), 
gave a
counterexample, and also stated what we call his
``conjecture''.

Hermann Weyl had given the modern definition of abstract 
differentiable
manifolds in about 1913 in his book on Riemann surfaces, 
yet Elie Cartan in his
groundbreaking 1925 book on Riemannian Spaces (today's 
Riemannian manifolds)
not only wrote, ``It is very difficult to define a 
Riemannian space,'' but in
fact never did define them.   

Now this book is a most marvelous
piece of mathematics which is full of unexplained and 
nonrigorous (for
most of us) terms such as ``infinitesimal rotation''.  Was 
this
speculative mathematics?  Was it criticized by 
contemporaries for its
lack of rigor? Should it not have been published?  In fact 
it 
is---now---perfectly easy to rigorize it using the theory of
connections in fibre spaces, invented after the book was 
published.
 
\lethead{More Cautionary Tales}
There are other Cautionary Tales to be told:  

\lethead{Even Gauss Published Incomplete Proofs}
Gauss's first proof of the Fundamental Theorem of Algebra, 
in his 1799
dissertation, was widely admired as the first wholly 
satisfactory
proof.\ It relied, however, on a statement ``known from 
higher
geometry'', which ``seems to be sufficiently well 
demonstrated'': {\it
If a branch of a real polynomial curve $F(x,y)=0$ enters a 
plane
region, it must leave it again}.  Gauss, evidently feeling 
more
persuasion was needed, added: ``Nobody, to my knowledge, 
has ever
doubted it. But if anybody desires it, then on another 
occasion I
intend to give a demonstration which will leave no doubt 
$\ldots$\!.''
According to Smale's 1981 {\it Bulletin\/} article (from 
which these quotes
are taken), this ``immense gap'' remained even when Gauss 
redid this
proof 50 years later, and the gap was not filled until 
1920.  

\paragraph{Simple Finite Groups}%
The authors allude to the 15,000 published pages 
comprising the
classification of finite groups as ``theoretical'' 
mathematics, and an
example of ``big science'' in mathematics, but they do not
characterize it as a Cautionary Tale, as I do.  Has the 
classification
been rigorously proved?  What kind of a proof is this?  Is 
there an
expert who claims to have read it all and verified it?  It 
is
overwhelmingly probable that 15,000 pages contain 
mistakes.  I have
been told recently that some of the parts that were farmed 
out by the
organizers of the project have never in fact been 
completed.  What
then is the status of the classification theorem?  Can we 
rely on it?
If in fact the proof is incomplete, shouldn't this be made 
public?
Who's in charge here, anyway?

\paragraph{Computer-assisted Proofs}%
These present many Cautionary Tales.  Oscar Lanford
pointed out that in order to justify a computer 
calculation as part of
a proof (as he did in the first proof of the Feigenbaum 
cascade
conjecture), you must not only prove that the program is 
correct (and
how often is this done?), but you must understand how the 
computer
rounds numbers, and how the operating system functions, 
including how the
time-sharing system works.  In fact, Lanford pointed out, 
my late
colleague R. Devogelaere discovered an error in Berkeley's 
system
caused by the time-sharing protocol.

\paragraph{The 4-color Theorem}%
A case in point, combining features of the two preceding 
Cautionary
Tales, is the proof by Appel and Haken of the 4-color 
theorem.  In
their interesting 1986 article in the {\it Mathematical 
Intelligencer\/} they point
out that the reader of their 1977 articles must face ``50 
pages
containing text and diagrams, 85 pages filled with almost 
2500
additional diagrams, and 400 microfiche pages that contain 
further
diagrams and thousands of individual verifications of 
claims made in
the 24 lemmas in the main section of the text.''  In 
addition the reader
is told that ``certain facts have been verified with the 
use of about
twelve hundred hours of computer time$\ldots$in some 
places there
were typographical and copying errors.'' They go on to 
assert that
readers of the 1986 article will understand ``why the type 
of errors
that crop up in the details do not affect the robustness 
of the
proof.''  They point out that every error found subsequent 
to publication
was ``repaired within two weeks''.  Several new errors and 
their
corrections are discussed in this article, as well as an 
``error
correction routine'' that seems to the authors ``quite 
plausible''.
In 1981 ``about 40 percent''  of 400  key pages had been
independently checked, and 15 errors corrected, by  U. 
Schmidt.  In
1984 S. Saeki  found\ another error which ``required a 
small change in
one diagram and in the corresponding checklist'' (p.\ 58). 
 Are we now to
consider the 4-color theorem as proved in the same sense 
as, say, the
prime number theorem?  

\paragraph{Computerized Mathematics}%
Large-scale programs such as Mathematica and Maple are 
increasingly
relied on for both numerical and symbolical calculations.  
A recent
widely distributed email message reproduced Maple output 
which, if
correct, disproves the unique prime factorization theorem. 
 What is
the status of a proof that includes such calculations?  
Are we to
consider any proof relying on Mathematica, Maple or other 
such
programs merely speculative mathematics?  Should such a 
proof be so
labeled?

\paragraph{Don't Prove, Just Lecture!}%
A great deal of time has been wasted by respected 
mathematicians who announce
the solution of a famous conjecture, sometimes with a 
great deal of publicity
in the popular press, and then lecture widely about it, 
but never giving
details of a proof either in lectures or print, and who 
eventually admit they
are wrong. (As a variant, a proof is published in an 
unrefereed article.)  This
is speculative mathematics at its worst, and is 
inexcusable.\footnote"$^1$"{A
peculiar converse phenomenon is that of a reliable 
mathematician raising doubts
about a long-accepted proof, then working hard (or 
assigning students) to
correct it, only to eventually find that in fact the 
original proof is correct.
I have seen this happen twice with Birkhoff's proof of 
Poincar\'{e}'s Annulus
Twist theorem.}

\lethead{Dated, Labeled Proofs}
Perhaps published mathematics, like good wine, should 
carry a date. If
after ten years no errors have been found, the theorem 
will be
generally accepted.  But there should also be an 
expiration time: If
any thirty-year period elapses without publication of an 
independent
proof, belief in the theorem's correctness will be 
accordingly
diminished.  This would allow for the reality that 
concepts of proof
and standards of rigor change.  

In addition we could attach a label to each proof, e.~g.:
computer-aided, mass collaboration, formal, informal, 
constructive,
fuzzy, etc.  Each theorem could then be assigned a number 
between zero
and one characterizing its validity, to be calculated from 
the proof's
label, vintage (see above), and the validities of the 
theorems on
which it is based.  These controversial calculations would 
themselves
become the subject of a new field of research$\ldots$

\lethead{Research Announcements}
I agree with the authors' list of problems with speculative
mathematics.  I think their prescriptions are sensible, 
{\it except}
for doing away with research announcements.  Research 
Announcements,
as published in this {\it Bulletin},
are a Good Thing!  In them I often read interesting 
accounts of
striking new results in fields I'm not expert in. It would 
not occur
to the authors to send me their results in preprint or 
electronic
form, since they don't know I'd be interested.  And I am 
sure that I
will not read the complete proofs---it is the clear and 
succinct
statement of important results, along with some indication 
of method,
{\it and some independent reassurance of correctness}, 
that I value in R.A.s.

The advantage of an R.A. over the complete proof is that 
it is
published more quickly.  A second advantage is that it is 
subject (in
current practice) to rigorous scrutiny by the editors, who I
understand insist on seeing some writeup of the proof.  
This can only
improve the complete version.  This of course is a lot of 
work for the
editors, and no doubt expensive for the AMS---but in my 
view R.A.s make
a unique and valuable contribution.

If we had {\it more} R.A.s there
would be less excuse for premature presentations of 
speculative
results---the rejoinder  would be, ``Why don't you publish a
Research Announcement?''

\lethead{References}

K. Appel and W. Haken, {\it Every planar map is four 
colorable.  Part \RM{I:}
Discharging}, Illinois J. Math. {\bf 21} (1977), 421-490; 
{\it Part
\RM{II:} Reducibility}, ibid, 491-567.

K. Appel and W. Haken, {\it The four color proof 
suffices}, Mathematical
Intelligencer {\bf 8} (1986), 10-20.

S. Smale, {\it The fundamental theorem of algebra and 
complexity
theory}, Bull. Amer. Math. Soc. {\bf 4} (1981), 1-36.

P. Heegard, {\it Sur l'analysis situs}, Bull. Soc. Math. 
France {\bf 44}
 (1916), 161-142. [Translation of: {\it Forstudier til en
topologisk teori f\"{o}r de algebraiske Sammenhang}, 
dissertation,
Copenhagen, Det Nordiske Forlag Ernst Bojesen, 1898.] 

H. Poincar\'{e}, {\it Sur les nombres de Betti}, Comptes 
Rendus Acad. Sci.
(Paris) {\bf 128} (1899), 629-630.  [Reprinted in \OE 
uvres {\bf 6}.] 

H. Poincar\'{e},   {\it  Compl\'{e}ment \`a l'analysis 
situs},
  	Rendiconti Circolo Matematico Palermo {\bf 13} (1899),
  285-343.  [Reprinted in \OE uvres {\bf 6}.]

H. Poincar\'{e}, {\it Deuxi\`{e}me compl\'{e}ment \`a 
l'analysis situs},
Proc. London. Math. Soc. {\bf 32} (1900), 277-308.  
[Reprinted in
\OE uvres {\bf 6}.] 

H. Poincar\'{e}, {\it Cinqui\`{e}me compl\'{e}ment \`a 
l'analysis
situs}, Rendiconti Circolo Matematico Palermo {\bf 18} 
(1904),
45-110.  [Reprinted in \OE uvres {\bf 6}.]


\newletter

\letaddress
Saunders Mac Lane
Department of Mathematics
University of Chicago
Chicago, IL  60637
saunders@MATH.UCHICAGO.EDU
\endletaddress

	In the July 1993 issue of this {\it Bulletin\/}, Arthur 
Jaffe and Frank Quinn
have speculated about a possible synthesis of mathematics 
and theoretical
physics.  On the way they make many interesting 
observations.  However, in my
view, their main proposal is both hazardous and 
misconceived.

	In the fall of l982, Riyadh, Saudi Arabia, was the seat 
of the First
International Conference on Mathematics in the Gulf 
States.  Michael Atiyah
attended, and provided outstanding advice to many of the 
younger conferees; I
admired his insights.  One evening, one of our local hosts 
gave an excellent
dinner for a number of the guests;  in his apartment we 
feasted on the best
lamb (prepared by his wife, who did not appear).  After 
this repast, we all
mounted to the roof of the apartment house, to sit at ease 
in the starlight. 
Atiyah and Mac Lane fell into a discussion, suited for the 
occasion, about how
mathematical research is done. For Mac Lane it meant 
getting and understanding
the needed definitions, working with them to see what 
could be calculated and
what might be true, to finally come up with new 
``structure'' theorems.  For
Atiyah, it meant thinking  hard about a somewhat vague and 
uncertain situation,
trying to guess what might be found out, and only then 
finally reaching
definitions and the definitive theorems and proofs. This 
story indicates that
the ways of doing mathematics can vary sharply, as in this 
case between the
fields of algebra and geometry, while at the end there was 
full agreement on
the final goal: theorems with proofs.  Thus differently 
oriented
mathematicians have sharply different ways of thought, but 
also common
standards as to the result.

	Jaffe and Quinn misappropriate the word ``Theoretical'' 
as a label for
what is really speculation.  This will not do; the word 
``Theory'' has a firm
mathematical use, as in the Theory of a Complex Variable 
or the theory of
groups.  In Jaffe-Quinn, I note also their complete 
misunderstanding of the
recently accomplished classification of all finite simple 
groups.  It was not
primarily a matter of organization or of ``program'', but 
one of inspiration,
from the early ideas of Richard Brauer (Int. Congress 
l954), the work of
Hall-Higman, the famous odd-order paper of Feit and 
Thompson, and the striking
discovery of new sporadic simple groups, as with the Janko 
group and the
Fischer-Greiss Monster, plus a decisive summer conference 
in 1976.  A small
part of the classification is not yet published---that for 
certain thin
subgroups. Here and throughout mathematics, inspiration, 
insight, and the hard
work of completing proofs are all necessary. No guide from 
physics can help,
and the occasional suggestion (Atiyah) that this 
classification should be all
done conceptually or geometrically has not (yet?) worked 
out.

	The sequence for the understanding of mathematics may be: 

\centerline{\bf intuition, trial, error, speculation, 
conjecture, proof.}

\noindent
The mixture and the sequence of these events differ widely 
in different
domains, but there is general agreement that the end 
product is rigorous
proof---which we know and can recognize, without the 
formal advice of the
logicians.  In many parts of geometry, differential 
topology, and global
analysis the intuitions are very complex and hard to 
reduce to paper; as a
result there can be a long development before closure.  An 
example is the
brief paper of Kirby establishing the annulus conjecture 
by using a sequence 
of known work to reduce the conjecture to a question of 
homotopy theory. In
each case, the ultimate aim is proof; for example, the 
review of a l988 paper 
by Goretsky and MacPherson in {\it Mathematical Reviews\/} 
states ``This long article
completes proofs of results that the authors have been 
announcing since l980.''
I surmise that much of the delay was needed to get matters 
in order, but the
old saying applies ``Better Late than Never'', while in 
this case ``never''
would have meant that it was not mathematics.  For the 
Italian geometers at 
the turn of the century the better late was very much 
later; the Italian
intuitions needed---and encouraged---the working out of 
many rigorous 
algebraic and topological methods by a long array of 
experts: van der  
Waerden, Krull, Zariski, Chevalley, Serre, Grothendieck 
and many others. 
Intuition is glorious, but the heaven of mathematics 
requires much more.    
As a result, Zariski and Grothendieck clearly outrank all 
the Italians. 
Mathematics requires both intuitive work (e.g., Gromov, 
Thurston) and 
precision (J. Frank Adams, J.-P Serre).  In theological 
terms, we are not 
saved by faith alone, but by faith and works.

	{\it Conjecture\/} has long been accepted and honored in 
mathematics,
but the customs are clear. If a mathematician has really 
studied the subject
and made advances therein, then he is entitled to 
formulate an insight as a
conjecture, which usually has the form of a specific 
proposed theorem.
Riemann, Poincar\'e, Hilbert, Mordell, Bieberbach, and 
many others have made
deep such conjectures.  But the next step must be proof 
and not more
speculation. On the Poincar\'e hypotheses, Henry Whitehead 
published an
erroneous proof but soon recognized the gaps; later, 
listening to proposed
proofs by others, he could say ``And now you do 
this\!\dots\!and then 
that\!\dots\<\<''
to the needed effect.  Sadly he was not on hand at the 
recent false claim made
for the Poincar\'e hypothesis in the {\it New York 
Times\/}, but other experts in
Kirby's seminar were at hand.  False and advertised claims 
have negative value,
even in these days of undue pressure to publish. The {\it 
New York Times\/}, in this
and other recent flamboyant cases, does not classify as a 
refereed journal.

	{\it Speculation\/}, unlike conjecture, usually is a much 
less specific
formulation  of some guess or insight; sometimes 
speculations can be combined
into a  program or outline of possible further work. In 
the l930s there was
such a program for the arithmetization of the class field 
theory, while the
current Langlands proposals provide such a program in 
non-abelian class field
theory and its connections to representation theory. 
Programs, good and
diffuse, come in all sizes and with very different 
prospects. Good programs
depend on insight; their execution requires proof.

	{\it Errors\/}, alas, abound; no one is immune, but the 
more egregious
ones often arise from overly hasty confidence that some 
insight can be filled 
in later by some technique.  On occasion, some leading 
journals ({\it The Annals 
of  Mathematics\/}) have been sometimes careless in not 
checking for errors in  
papers on fashionable subjects.  Moreover, when error is 
discovered, the
journal should publish a retraction, so that all workers 
may know.  Some
notable errors, such as Dehn's lemma, have been a stimulus 
to subsequent work,
although Max Dehn himself, in this country after l933, did 
not profit.  In  
one more recent case a published error was overcome in a 
long paper of  
several hundred pages.  The reviewer remarked ``This paper 
is written in  
great detail, sometimes almost too much, but this is not a 
bad thing, given 
the long record of incorrect proofs in this subject.'' 
This states clearly 
the goal: It is not mathematics until it is finally proved.

	The recent fruitful interchange of ideas (connections, 
fiber bundles,
etc.)  with physics (quantum gravity and all that) has 
been a decided stimulus
and a source of new ideas and reapplication of old ones.  
It's great, but
involves some of the current weaknesses of physics. Thus, 
when I attended a
conference to understand the use of a small result of 
mine, I heard lectures
about ``topological quantum field theory'', without the 
slightest whiff of a
definition; I was told that the notion had cropped up at 
some prior conference,
so that ``Everybody knew it.'' Much the same may apply to 
``Quantum groups'',
which are not groups.  This practice reflects 
carelessness, sloppy thinking,
and inattention to established terminology, traits which 
we do not need to copy
from the physics community.  There, it was said that one 
person asked another
``What are you working on?''; ``String theory''; ``Oh, 
didn't you know, that
went out of fashion last week.'' We are fortunate that 
mathematics has a more
permanent character and is not (at least yet) bound to 
concentrate everybody in
the latest thing. The Lord's house has many mansions.

	For other and deeper reasons I cannot share the 
enthusiasm of
Jaffe-Quinn for physics. Their comparison of proofs in 
mathematics with
experiments in physics is clearly faulty. Experiments may 
check up on a theory,
but they may not be final; they depend on instrumentation, 
and they may even be
fudged.  The proof that there are infinitely many 
primes---and also in suitable
infinite progressions---is always there. We need not sell 
mathematics short,
not even to please the ghost of Feynmann.

	Since World War II, physics has played a dominant role in 
American
science.  But today, it faces serious troubles. The 
standard theory requires
the existence of a ``Higgs Boson'', which has not yet been 
found; searching for
it requires the Superconducting Supercollider, costing 
billions and requiring
annual appropriation.  Then Relativity Theory has led to 
plans for tests of the
existence of gravity waves in an expensive LIGO (Laser 
Interferometry Gravity
Wave Observatory). The study of the Big Bang appears to 
mix speculation and
science. Senior physicists have time to write popular 
books on the ``Final
Theory''. For younger physicists it can be hard, with 
problems depending on the
funding for big apparatus or on papers with l90 authors.  
One such youngster
(David Lindler) recently broke ranks to write a book {\it 
The End of Physics\/};
his main contention is that cosmologists and theoretical 
physicists encourage
each other to wider and wilder speculations.

	  Physics  has provided mathematics with many fine 
suggestions and new
initiatives, but mathematics does not need to copy the 
style of experimental
physics. Mathematics rests on proof---and proof is eternal.


\newletter

\letaddress
Benoit B. Mandelbrot
Mathematics Department
Yale University
New Haven, CT 06520-8283
fractal@watson.ibm.com
\endletaddress

  Unfortunately, hard times sharpen hard feelings; witness 
the
discussion (to be referred to as JQ) that Arthur Jaffe and 
Frank
Quinn have  devoted to diverse tribal and territorial issues
 that readers of this {\it Bulletin} usually leave to 
private gatherings.
Those readers --- you! --- like to be called simply 
``mathematicians''.  
But this term will not do here, because my comment is 
ultimately founded on
the following conviction:

{\narrower\smallskip
For its own good and that of the sciences, it is critical 
that
mathematics should belong to no self-appointed group;  no 
one has, or
should pretend to, the authority of ruling its 
use.\smallskip}

  Therefore, my comment needs a focussed term to denote 
the typical
members of the AMS. Since it is headquartered on Charles 
Street, I 
propose (in this comment, and  never again) to use 
``Charles mathematicians''.

The reason I agree to respond to JQ is that informal 
soundings
suggest that---save for minor reservations---Charles 
mathematicians
tend to admire it.  To the contrary, I find most of it 
appalling.
Least among my reasons is the manner in which JQ refers to 
my work (to
this, I shall respond toward the end of this piece).

My main objection to JQ is that, in their search for 
credit for some
individuals at the expense of others they consider rogues, 
they propose 
to set up a police state within Charles mathematics, and a 
world cop
beyond it borders.  How far beyond?  They do not make it 
clear.
Indeed, ``technical mathematics'' is not only a silly 
term, but also empty,
 unless it is applied to nearly every field outside Charles
mathematics.  Everyone is enchanted when seemingly
disparate problems turn out to have a common solution.  How
attractive, therefore, if the same treatment could be 
applied to all
``theoretical mathematicians''!  But politics is more 
complicated
than any science, and the term ``theoretical 
mathematicians'' 
does not apply to any self-defined entity. Therefore, this 
discussion of 
JQ must deal separately with Charles mathematicians, and 
with others.

Let us begin with the latter. In their concern about the
relations between mathematicians and physicists, JQ 
bemoans that
``students in physics are generally indoctrinated with
anti-mathematical notions and\dots often deny their work is
incomplete.''  To change the situation, JQ does not 
propose anything
like a negotiation between equals, only a prescription 
(their word).
They tell those foreigners which proper behavior would 
increase the
happiness of a few Charles mathematicians.  Clearly, the 
prescription
will only be heard by the very few already conditioned to 
heed it.
Why should scientists care about credit to Charles 
mathematicians, 
given the Charles mathematicians' own atrocious record in 
giving credit.  
The pattern---sad to say---is an acknowledgment that 
``it is natural to take this tack or that'' (when a 
physicist has 
sweated to establish that this is indeed ``natural'') or 
an acknowledgment 
that ``physicists believe'' (which physicists?) or ``the 
computer has shown'' 
(the computer by itself? or perhaps some unmentioned 
``technician'').

However, the main reason why I find the JQ prescription 
appalling is
because it would bring havoc into living branches of 
science.  Philip Anderson
describes mathematical rigor as ``irrelevant and 
impossible''.  I would
soften the blow by calling it ``besides the point and 
usually distracting,
even where possible''.  As a first example, take the 
statistical process of
percolation.  It involves a clearly mathematical 
construction discovered
by Hammersley, a mathematician, but it soon fell into the 
hands of
theoretical physicists.  They gradually discovered that 
percolation
illustrates hosts of natural
phenomena of interest to them, and they established an 
extraordinarily
long list of properties (many of which, may I add, are 
fractal in
character).   In the meantime, the mathematicians (and 
their numbers
include truly distinguished persons) have lagged way 
behind.  They have 
brought little that my physicist friends find valuable.  
Anyone who will give rigorous solutions to this enormously 
difficult
problem will deservedly be hailed in mathematics, even if 
the work
only confirms the physicists' intuition.  I hope that 
proper credit
will be given to individual physicists  for
their insights, but fear it will not.  Will
this mathematics be also noticed by phycisists?  Only if 
they prove
more than was already known, or if the rigorous proofs are 
shorter
and/or more perspicuous than the heuristics.

Needless to say, everyone knows that  mathematical physics 
overlaps
with a community of physicists devoted to ``exact 
results''.  But,
for all practical purposes, these brilliant people have 
crossed over
to become Charles mathematicians.  
The same has happened to whole disciplines.  Take 
mathematical
statistics.   Jerzy Neyman disciplined his followers to 
practice the
most exacting 
rigor (even though their results are of interest to hardly 
any
Charles mathematician).  Now that the long shadow of 
Neyman has
waned, the mood has changed, and mathematical statistics 
has been
freed to seek a place in the community of sciences.  But 
if the JQ 
prescription were applied to it, it would die once again.

Now let me move back from foreign to internal affairs.  
Many Charles
mathematicians are offended that a few major players have 
been well
rewarded for ``passing'' statements as if they were
proven theorems.  The fear is expressed, that the person 
who will
actually prove these  assertions will be deprived of 
credit.  I
think this is an empty fear because the AMS already has 
countless
ways of rewarding or shunning whomever it chooses.  Even 
if one could
agree with the goals of JQ, I would contend that, in order
to avoid ``problems'' caused by current custom, there is 
the need of
altering a community's present rules of behavior.
  Moreover, I do not agree with the goals of JQ, because I 
find them
destructive.  My reading of history is that mankind 
continually
produces some individuals with the highest mathematical 
gifts who
will not (or cannot) bend to pressures like those proposed 
by JQ.  If
really pressured, they will leave mathematics---to 
everyone's great loss.

My first witness will be the probabilist Paul L\'evy 
(1886-1971).  (Moving
back in time adds perspective!)  The French-style Charles
mathematicians of his time kept blaming him for failure to 
prove
anything fully (and for occasional lapses in elementary
calculations!).  There was no field where he could flee 
away from
Charles mathematics, but he did not change.   He went on, 
well
into his seventies, producing
 marvelous  and startling intuitive ``facts'' that may 
have been
``incomplete'' yet continue to provide well-rewarded work 
for many.
Yet, when he was 71 (and I was a junior professor working 
for him), he
continued to be prevented from teaching probability, part 
of a pattern of tormenting him in every way conceivable.   
Question:
Who gained?

My second witness will be Poincar\'e himself.  In the 
recently 
published letters from Hermite, his mentor,  to
Mittag-Leffler, there are constant complaints about 
Poincar\'e's
unwillingness to heed well-intentioned advice and polish 
and publish
full proofs.  Concluding that Poincar\'e was incurable, 
Hermite and E.
Picard (who inherited his mantle) shunned
Poincar\'e, prevented him from teaching mathematics, and 
made him
teach mathematical physics, then astronomy.  His published 
lecture
notes cover basic optics, thermodynamics, 
electromagnetism: a second-
or third-year ``course in theoretical physics''.  It may 
well be true
that Poincar\'e's 1895
{\it Analysis Situs} remained for quite a while a ``dead 
area''.
Questions:  Who was harmed?  Would the world have been a 
happier 
place if Poincar\'e, awed by Hermite and his cohorts, had
waited to publish until he knew how to irrigate this dead 
area?

One could bring other witnesses, but few as great as 
Poincar\'e and
L\'evy.  Why were there so few like them during the recent 
period?  One
possible reason lies in the flow of young people who kept
being introduced to me for advice.  They were acknowledged 
as brilliant
and highly promising; but they could not stomach the 
Bourbaki
credo, hence saw no future for themselves in mathematics.  
Bystanders
who understood what was happening asked me to help these 
young people
to hang on, but there was no way.  They did not straighten 
up but scattered, 
at a loss, as I see it, both for themselves and for 
mathematics.

  Before concluding this response to JQ, let me say that I 
am
disappointed that JQ should mention me, without adding that 
whenever I see in my work something that might interest 
Charles 
mathematicians, I make it a point of seeking out their 
attention and
describing what I had done as a conjecture, therefore a 
challenge to be
either proven or disproven, within the usual standards of 
rigor.
It is a delight that---in all the branches I touched, 
harmonic analysis,
probability theory and (most widely known) the theory of 
iteration of
functions---brilliant people promptly took up my 
challenge, and were
led to beautiful theorems. 

   What about credit, an issue that dominates JQ?  For the 
proofs
of my conjectures, full credit is due, and no one denies it.
For pioneering the use of computer graphics in mathematics,
raising the problems and making early conjectures, 
full credit is also due.  But to my regret
many Charles mathematicians extend it with undisguised 
reluctance.
Coming from my uncle Szolem (1899-1993) and other persons 
I like to
admire without reservation, this reluctance used to be 
annoying.  
Luckily, I have long known that it does me no harm, thanks 
to public 
acclaim from sources free of the biases and hangups of 
Charles mathematics.

   For this reason, I am dismayed that JQ should bring in 
S. C. Krantz as
``an expression of the mathematical discomfort with my 
activity.''  Beyond  ``expressing discomfort'',
Krantz also tries to justify it by  anecdotes (one of his
specialties) based on imagined 
 ``facts'' and wrong dates, as I have shown in the 
foreword I wrote to the book \it Fractals for the Classroom,
Part I\rm, by H-O. Peitgen, H. J\"urgens, and D. Saupe 
(Springer 1991). Had
I not expected Krantz's press releases and other utterings 
to sink
promptly into oblivion, I would have answered in a more 
public forum.
But, since his opinion has surfaced in JQ it is important to
 refer to my full rebuttal. 

   In its concentration on credit and its overbearing 
attitude towards proven 
and successful domains that are clearly outside of Charles 
mathematics, 
the JQ piece is in dreadful bad taste.  
Its complaint that rigor is threatened is no more true today
than it was in the midst of Bourbaki's domination.  Its 
attitude towards
Charles mathematician rogues would, if it returns to 
power, deprive
us of future Poincar\'es and Levys.  Let me, however, end 
on a
positive note.  It is a pleasure that JQ do not put into 
question the
recent improvements in the mood between mathematics and its
neighbors.  The best borders are open borders that allow 
nominal
physicists to be praised for their mathematics and nominal
mathematicians to be praised for their physics.


\newletter

\letaddress
David Ruelle
Institute des Hautes Etudes Scientifiques
35 Rue de Chartres
91440 Bures-sur-Yvette
France
\endletaddress

Dear Dick,\par

   Thank you for soliciting my opinion on the Jaffe-Quinn 
paper.  I am glad
that this paper will appear in the BAMS, because it raises 
issues that deserve
to be discussed.  Since my own feelings are not extremely 
strong however, let
me just express a few comments in the form of this letter 
to you, from which
you will feel free to quote or not to quote.

   Nobody will question the need to indicate if something 
presented as a
``proof'' is really a complete, rigorous, mathematical 
proof, or something else
(to be specified).  For example in the study of the 
``Feigenbaum fixed point''
there is no quarrel: Feigenbaum's argument is profound and 
convincing, but
does not claim to be mathematically rigorous.  In other 
cases things are
less clear.  Kolmogorov did not publish a proof for the 
``K'' of ``KAM'', but
Yasha Sinai tells me that Kolmogorov gave lectures 
amounting to a full proof of
``K'' in the case of two degrees of freedom.

   One point that perhaps deserves being stressed is the 
usefulness of
cultural cross-fertilization in mathematics.  Feigenbaum's 
cultural background
in theoretical physics has allowed him to discover a new 
generic bifurcation
of smooth dynamical systems, which would not have been 
encountered soon by
following standard mathematical paths.  Similarly, the 
physical ideas of
equilibrium statistical mechanics have richly contributed 
to the mathematical
theory of smooth dynamical systems (with the concepts of 
entropy, Gibbs states,
etc., see my note in BAMS(NS) {\bf19} (1988), 259-268.  
The importance for pure
mathematics of ideas coming from theoretical physics is of 
course well known
to Jaffe and Quinn, and they are right in insisting that, 
with a little bit
of care, mathematics can benefit from these ideas without 
paying an exorbitant
price.


\newletter

\letaddress
Albert Schwarz
Department of Physics
University of California at Davis
Davis, CA 95616
asschwarz@ucdavis.edu
\endletaddress

   I agree completely with A. Jaffe and F. Quinn that 
heuristic
(``theoretical") work can be very useful for mathematics, 
but that it is
necessary to establish the rules of interaction between 
heuristic and rigorous
mathematics.  However I would like to suggest somewhat 
different rules. I'll
begin with a short exposition of my opinion and give some 
explanations at the
end of my letter as footnotes.

  1. I don't think that the name ``theoretical 
mathematics" as a common name
for heuristic mathematics and theoretical physics is 
appropriate. A. Jaffe and
F. Quinn consider only the part of theoretical physics 
that is not closely
related to the experiment, but this does not change the 
picture. The main
aim of theoretical physics is to explain experimental data 
or to predict the
results of  experiments. Many physicists believe now that 
string
theory can explain all interactions existing in nature. 
However today they are
not able to extract reliable predictions from  string 
theory because this is
connected with enormous mathematical difficulties.  The 
physicists have chosen
the only possible way: to analyze carefully the 
mathematical structure of 
string theory. This approach led to beautiful mathematical 
results, but still
did not give the possibility to solve realistic 
experimental problems. (*)

  String theorists believe that finally they will be able 
to give a formulation
of TOE (``theory of everything") on the base of string 
theory and to make all
calculations on this base. Other theoretical physicists 
have a doubt that
string theory is related to the physics of elementary 
particles at all.
However in any case the string theorists and all other 
theoretical physicists
have the same goal and the same psychology. The fact that  
the string theorists
use complicated mathematical tools and that their papers 
contain very important
contributions to pure mathematics does not mean that they 
became
mathematicians. Of course at a personal level the 
borderline is not so
sharp. There are scholars that can be considered as 
physicists and as
mathematicians at the same time; there are many papers 
that belong
simultaneously to theoretical physics and heuristic 
mathematics.

   Mathematics is not necessarily characterized by 
rigorous proofs. Many
examples of heuristic papers written by prominent 
mathematicians are given
in [1]; one can list many more papers of this kind. All 
these papers are
dealing with mathematical objects that have a rigorous 
definition.(**) 
However a mathematician reading a textbook or a paper 
written
by a physicist discovers often that the definitions are 
changed in
the process of calculation. (This is true for example for 
the definition of
scattering matrix in quantum field theory.)

   It would be meaningless to try to establish formal rules
of interaction between mathematics and theoretical physics.
Therefore I'll talk only about rules of interaction 
between rigorous
and heuristic mathematics. This means that I always have 
in mind that 
the papers under consideration are based on rigorous 
definitions.

   2. It is suggested in [1] that the ``theorists" should 
label their 
statements as conjectures. I agree that the word 
``theorem" should be reserved
for rigorously proven statement, but it would be arrogant  
to insist that
heuristic mathematics can produce only conjectures. The 
``conjectures"
in the sense of [1] are of very different nature in the 
range from wrong to
completely reliable
 (***). I believe that in heuristic papers one  should 
avoid the word
``theorem", using instead the words ``pretheorem", 
``fact", ``statement" and, 
of course, ``conjecture". I borrow the word ``pretheorem" 
from the book [2]
together with the explanation: The word ``pretheorem" 
should mean that there
exists at least one set of technical details supporting a 
theorem of the sort
sketched. The word ``fact" could mean statement that is 
not proven rigorously,
but is considered by the author as reliable (****). The 
word ``statement" can
denote an assertion provided with a heuristic proof that 
cannot be extended to 
a rigorous proof without essential new ideas. Finally, the 
word ``conjecture"
should be understood as a statement supported not by a 
heuristic proof, but
by an analogy, examples, etc.

   The suggestion of [1] about flags indicating
 heuristic character of the paper sounds reasonable. 
However, as I
mentioned already, heuristic papers are of a different 
nature, therefore
these flags should be of different colors.(*****)
This would be in complete agreement with the Jaffe-Quinn 
suggestion, but
more acceptable for authors of heuristic papers.

  3. There is no doubt that a mathematician that gave a 
rigorous proof of a
statement heuristically proved by another mathematician or 
physicist deserves
essential credit. However  I don't think that one can say 
a priori that ``a
 major share of credit for the final result must be 
reserved for the
 rigorous work."  Sometimes this
is true. For example in constructive quantum field theory 
rigorous proofs
are often a hundred times harder than heuristic 
considerations. In this case
complaints that physicists underestimate the work of 
mathematicians are
completely justified. However sometimes the final stroke 
requires only careful
student's work. (By the way, filling gaps in a heuristic 
paper of this kind or
in a research announcement would be an extremely valuable 
part of a student's
education.)

 4. It is important to stress that heuristic mathematics 
is a legitimate
part of mathematics. ( Now this is recognized only in 
mathematical physics.
A crucial role in this recognition was played by the 
policy of the
editorial board of {\it Communications in Mathematical 
Physics\/} headed by
A. Jaffe.) In the ideal case every rigorous mathematician 
has to work
also heuristically and every author of a heuristic paper 
fully recognizing the
importance of rigor has to try to prove his results 
rigorously. Every
mathematician begins his work with heuristic 
considerations, as emphasized in
[1]. If he obtained an interesting result, he should try 
to give a rigorous
proof of it. However if he did not succeed in this 
attempt, he can publish the
result with a heuristic proof with a hope that somebody 
else will be able to
finish the work. In our age of division of labor it is 
difficult to understand
why both parts of the work must be necessarily done by the 
same person.
Nevertheless most mathematicians are tied with an odd bias 
that
only rigorous results deserve publication. I hope that the 
discussion initiated
by Jaffe and Quinn will help to destroy this bias and 
instead of separation we
will see a community of scholars united by a common goal 
and sometimes acting as
rigorous mathematicians (if possible), sometimes writing 
heuristic papers (if
rigorous methods do not work).

   (*) Such a situation is not completely new in physics. 
It took about
 20 years to give a correct formulation of mathematical
 problems arising in quantum electrodynamics and to solve 
these problems in
the framework of perturbation theory. (I have in mind the 
invention of
renormalization theory.) The way from the introduction of 
gauge fields to the
quantization of these fields and to the construction of 
gauge theories giving 
a description of electromagnetic, weak and strong 
interactions was only a 
little bit shorter.

   (**) This does not mean that all definitions are 
formulated precisely.
Sometimes it is convenient to leave some freedom in a 
definition used in a
heuristic paper.  However this means only that we work 
with several rigorous
definitions at the same time.

   (***) Let me give an example. In 1987 I conjectured 
that the Jones polynomial 
can be obtained from quantum field theory. This conjecture 
was inspired by the
conversation with V. Turaev. Turaev told me that V. Jones 
invented an invariant
of knots, that can be considered as a generalization of 
the Alexander polynomial.
He thought that the application of the general method of 
construction of
topological invariants by means of quantum field theory 
suggested by me in
1978 [3] could give an explanation of the origin of the 
Jones invariant. (The
Alexander polynomial can be expressed in terms of 
Reidemeister torsion and
my paper contained construction of the smooth version of
the  Reidemeister torsion, the so-called Ray-Singer 
torsion.)
 Answering Turaev's question, I found a Lagrangian
(Chern-Simons Lagrangian) giving invariants of 
three-dimensional manifolds
and conjectured that these invariants are connected with 
the Jones polynomial [4].
A year later a heuristic proof of this conjecture was
 given in a brilliant paper [5] by Witten (who knew my 
paper of 1978, but did
not know the paper [4]). Of course, Witten went much 
further than I. (I
consider the contribution of [4] as negligible in 
comparison with [5] or [3].)
The constructions of his paper (in particular, the 
connection with
two-dimensional conformal theory) were exploited later in 
hundreds of papers
and led not only to heuristic, but also to important 
rigorous results. But in
the terminology advocated in [1]\ the difference between 
Witten's and my
statements disappears: both are qualified as conjectures!

  (****) Mathematicians often underestimate the 
reliability of heuristic proof.
Probably, the results of a good mathematician, working 
heuristically, are not
less reliable than the results of an average rigorous 
mathematician. (One can
 consider this statement as a definition of a good 
mathematician.)
 Mathematicians know that a formal proof leads always to a 
correct result, 
but they forget sometimes that they are human beings 
(``Errare
 humanum est"). Therefore erroneous ``rigorous" papers are 
not so rare. 
Heuristic methods are not completely reliable, therefore 
the scholars using
these methods have to apply all possible checks to 
guarantee reliability of
their results. However as stressed in [1] a rigorous proof 
gives often  new
insights and new results, therefore it is necessary also 
in the case when the
statement is completely reliable. Let me illustrate this 
fact by the following
example. Gauss found by direct calculation that the length 
of a certain curve
coincides with great accuracy with some arithmo-geometric 
average. Of course,
he conjectured that these two numbers coincide precisely. 
Calculating more and
more digits, he could make the probability that this 
conjecture is violated as
low as wanted. Instead he gave a rigorous proof 
discovering some properties of
elliptic functions, that definitely are much more 
interesting than the
conjecture itself. Probably, we will have a similar 
situation with some
conjectures about mirror symmetry. The coincidence of some 
numbers predicted by
mirror symmetry was checked in so many cases, that it is 
almost impossible to
doubt the correctness of these conjectures . Nevertheless 
neither physicists
nor mathematicians are satisfied. They would like to know 
the reason for this
coincidence.

 (*****) I would like to list some of the possible cliches.

$A$. The paper contains a complete rigorous proof.

$B$. The paper contains no proofs (or: The proofs are only
sketched), but the author gave detailed rigorous proofs of 
all results of the
paper.

 $B_1=B+$ The proofs  are written and available upon 
request.

 $B_2=B+$ The author is planning to write the proofs not 
later than...

 $B_3=B+$ The author is not planning to write detailed 
proofs. He
would be ready to help anybody willing to perform this work.

 $C$. The paper is addressed to physicists. Therefore the 
results are
not formulated as mathematical theorems. However it is 
easy to give
conditions making the proofs completely rigorous.

 $D$. The proofs in the paper are rigorous, but they are 
based on some
statements of the paper\dots having only heuristic proofs.

 $E$. The proof given in the paper cannot be considered as 
complete,
but the author believes that the gaps in the proof can be 
filled in
without essentially new ideas.

$E_i=E+\ldots$

 $F$. We give a heuristic proof of our statements. A 
rigorous proof
cannot be obtained by our methods.

A. Jaffe and F. Quinn propose to publish research
announcements only in
 the case when the complete paper is already written and 
refereed. They think
 that this is possible because the announcements can be 
distributed via e-mail
 as preprints. However one can use e-mail to distribute a 
complete paper
too! I
believe that pretty soon scientific journals will publish 
only research
 announcements (together with information on how to get 
complete proofs via
 e-mail), review papers and extremely important papers. 
For me personally it is
 easier already to find a paper that I am interested in by 
means of an
 electronic bulletin board, than in the library.

I am indebted to J. Hass, C. Tracy, and especially to Yu. 
Manin for very
useful discussions.

\lethead {References}

1. A. Jaffe, F. Quinn, {\it Theoretical Mathematics\RM{:}
 Toward a cultural synthesis
of mathematics and theoretical physics\/}, 
Bull. Amer. Math. Soc. {\bf 29}
(1993), 1-13.

2. J. Adams, {\it Infinite loop spaces\/}, Princeton, 1978.

3. A. Schwarz, {\it The partition function of degenerate 
quadratic functional
and Ray-Singer invariants\/}, Lett. Math. Phys. {\bf 2} 
(1978), 247-252.

4. A. Schwarz, {\it New topological invariants arising in 
the theory of
quantized fields\/}, Baku International Topological 
Conference,
Abstracts (Part 2) Baku, 1987.

5. E. Witten, {\it Quantum field theory and Jones 
polynomial\/}, Comm. Math.
Phys. {\bf 121} (1989), 351-399.


\newletter

\letaddress
Karen Uhlenbeck
Department of Mathematics
University of Texas at Austin
Austin, TX 78712-1082
uhlen@math.utexas.edu
\endletaddress

	The article by Arthur Jaffe and Frank Quinn has been a 
dynamic, healthy
catalyst for many interesting discussions about 
mathematics. I am very much of
the conviction that mathematics is much more than the bare 
and beautiful
structure as exposed by Bourbaki and as appreciated by 
myself before I had
research experience. Interest in mathematics from a 
broader-than-usual
perspective is presumed and advocated by the authors.

	I agree with many of the points of the article. Pure 
mathematicians
really ought to prove their theorems and publish their 
results in a clear and
understandable paper written in a timely fashion. What we 
may need, in addition
to the ``mathematica rejecta'' journal dreamed of in my 
youth, is the
``mathematica culpa'' elder journal? Some of the younger 
mathematicians could
be sent there for extreme sins, as well as us older folk 
who tend to end up
here as a way of life.

	My main criticism of the article is that it draws broad 
conclusions
from too narrow a perspective.  The relationship between 
physics and
mathematics has been fundamental to both for a long time. 
The gap between the
two is significant primarily in this century, as pure 
mathematics became very
abstract, experimental physics became very expensive, and 
the world became more
complicated.  However, even through this century, 
mathematics has relied on
physics for input quite steadily. I am sure other replies 
will point out the
influence of mechanics on calculus, optics on Riemannian 
and symplectic
geometry, general relativity on differential geometry, 
quantum mechanics on
functional analysis, geometric optics on harmonic
analysis,\ and gauge theory on
four-manifold topology.  This did not take place 
``neatly''.  The list would be
much longer if we included input from all sciences.

	Hence, I feel that the article makes exactly the wrong 
point about
influence on young mathematicians.  I well remember that 
as an undergraduate I
was initiated into the mysteries of distributions by being 
told by a graduate
student that physicists had used them, but understood 
nothing important about
them.  Only an innovative and brilliant mathematician like 
the idolized Laurent
Schwartz could make sense of the physicists' nonsense.  
Unfortunately, this
attitude was reinforced during my formative years by both 
mathematicians and
physicists.  Mathematicians seemed to think that 
physicists did not do physics
``right'', while physicists thought of mathematicians as 
worthless insects.
Only after taking part in the mathematical development of 
gauge theory could I
comprehend the essential importance of outside ideas in 
mathematics and the
contrary possibility of mathematical language being of 
real use outside the
discipline itself.

	I find it difficult to convince students---who are often 
attracted into
mathematics for the same abstract beauty and certainty 
that brought me
here---of the value of the messy, concrete, and specific 
point of view of
possibility and example.  In my opinion, more 
mathematicians stifle for lack of
breadth than are mortally stabbed by the opposing sword of 
rigor.

	As you can see, in the first part of my answer, I 
basically agree with
all the premises of the article.  I have serious 
objections of another sort to
the idea of creating a discipline called ``theoretical 
mathematics''.  Setting
aside the semantics, in the broader context of its 
description, ``theoretical
mathematics'' already exists.  It is called ``applied 
mathematics'', a much
bigger field than pure mathematics.  Applied mathematics 
is done mostly outside
departments of mathematics and draws in far more resources 
and many broad
scientific interests.  Only the combined elitism of very 
pure mathematics and
high-energy fundamental physics would claim that its own 
brand of speculative
and applicable mathematical structure should have a 
special name.  Would
nonlinear dynamics, which has an active and interesting 
interface with other
sorts of physics, qualify as theoretical mathematics?  
What about mathematical
biology, which may be held back by lack of mathematical 
attention to handling
complex information.  Some claim this field desperately 
needs mathematical
insight.  I would very much like to see the dialogue 
started by Jaffe and Quinn
extended to cover glories and disasters of interaction 
between pure mathematics
and the many other more applied areas of relevance.

	In conclusion, pure mathematicians might well spend even 
more time
building intellectual bridges to the rest of the 
scientific world. Jaffe and
Quinn imply that it would help to collect a toll for 
crossing one of few
well-built bridges.  They have, however, done a great 
service by describing it
in detail as worthy of tariff.
\vfill


\newletter

\letaddress
Ren\'e Thom
Institute des Hautes Etudes Scientifiques
35 Rue de Chartres
91440 Bures-sur-Yvette
France
\endletaddress

Dear Dick,\par

  Many thanks for your letter of May 21st with the 
enclosed article by Arthur
Jaffe and Frank Quinn.  I have many reasons to be 
interested in it, not only
because I am personally implicated in the ``Cautionary 
Tales''.  
There, I
can only confirm that the description of my evolution with 
respect to
mathematics is fairly accurate.  Before 1958 I lived in a 
mathematical milieu
involving essentially Bourbakist people, and even if I was 
not particularly
rigorous, these people---H.~Cartan, J.-P.~Serre, and 
H.~Whitney (a would-be
Bourbakist)---helped me to maintain a fairly acceptable 
level of rigor.  It was
only after the Fields medal (1958) that I gave way to my 
natural tendencies,
with the (eventually disastrous) results which followed.  
Moreover, a few years
after that, I became a colleague of Alexander Grothendieck 
at the IHES, a fact
which encouraged me to consider rigor as a very 
unnecessary quality in
mathematical thinking.  I somewhat regret that the 
authors, when quoting my
work in singularity theory, did not emphasize its positive 
aspects, namely, the
transversality lemma (with respect to jet systems), the 
theory of stratified
spaces (allowing for some anticipatory work by H.~Whitney 
and S.~\L
ojasiewicz), the characterization of ``gentle maps'' 
(those without blowing
up), the II and III isotopy lemmas.  All this was {\it 
written\/} for the first
time in my unrigorous papers.  Of course many people 
(Milnor, Mather,
Malgrange, Trotman and his school, McPherson, to quote 
just a few) may claim to
have a large part in the rigorous presentation of this 
theory.

        This leads me to the Jaffe-Quinn paper  itself, 
which involves a
very important question, and provides, I think, the first 
occasion (apart from
some solemn observations of S.~Mac Lane) for an in-depth 
discussion on
mathematical rigor.  I do still believe that rigor is a 
relative notion, not an
absolute one.  It depends on the background readers have 
and are expected to
use in their judgment.  Since the collapse of Hilbert's 
program and the advent
of G\"odel's theorem, we know that rigor can be no more 
than a local and
sociological criterion.  It is true that such practical 
criteria may frequently
be ``ordered'' according to abstract logical requirements, 
but it is by no
means certain that these sociological contexts can be {\it 
completely\/}
ordered, even asymptotically.

        One main argument of the Jaffe-Quinn paper is that 
we have to know,
when we want to use it for further research, if a 
published result may be
considered as ``firm'' as another, whether its validity 
may be universally
accepted.  My feeling is that it is unethical for a 
mathematical researcher to
use a result the proof of which he does not ``understand'' 
(except for the
specific case where he wants to disprove the result).  In 
principle, of course,
understanding here means a thorough knowledge of all the 
arguments involved in
the written proof.  From this viewpoint, it may not be as 
necessary as is
usually thought to classify all known truths in a 
universal library.  But
finally I think the proposal of the authors, to establish 
a ``label'' for
mathematical papers with regard to their rigor and 
completeness, is an
excellent idea.

        Rigor is a Latin word.  We think of {\it rigor 
mortis\/}, the rigidity
of a corpse.  I would classify the (would-be) mathematical 
papers under three
labels:\par

{\baselineskip=12pt\parindent=20pt\narrower

    \item{1)} 
\hbox{\hskip5pc\vbox{\vskip2.5pc}}
a crib, or baby's cradle, denoting
		``live mathematics'',
                allowing change, clarification, completing 
of proofs,
                objection, refutation.
\vskip 3pt
    \item{2)}  
\hbox{\hskip5pc\vbox{\vskip2.5pc}}
the tombstone cross.  Authors pretending 
		to full rigor,
                claiming eternal validity, may use this 
symbol as freely
		as they
                wish. This kind of work would constitute 
``graveyard
		mathematics''.
                
\vskip 3pt
   \item{3)}    
\hbox{\hskip5pc\vbox{\vskip3pc}}
the Temple. This would be a label delivered by
		an external
                authority, the ``body of high priests''.  
This body could
                initially be made up of the editors in 
chief of the ``core''
                papers as suggested by Jaffe-Quinn.  Its 
task would be to
                bestow the label at least on those
\crossymbol papers
                with sufficient promise to justify close 
examination. Later
                on, the IMU could decide on a permanent 
procedure to establish
                the priestly body, allowing for a 
relatively quick
                turnover of people in charge, with 
equitable worldwide 
                geographic representation.  One might 
suppose that such an
                institution could last a very long time.  
Should it however
                eventually come to grief, the unattainable 
nature of absolute
                rigor would be thereby demonstrated.\par}

        Let me end with a personal observation.  The 
Jaffe-Quinn paper
discusses at length the situation of mathematical physics, 
but does not seem to
admit that the problem may arise in other disciplines for 
which (unlike
physics) E. Wigner's phrase about the ``unreasonable 
effectiveness of
mathematics'' is not valid.  I strongly disagree with such 
a restriction.  I
see no reason why mathematics (even without computers and 
numerical
computation) should not be applied in other disciplines, 
in biology for
example.  In particular I believe that there are in 
analytic continuation
singular circumstances (unfoldings, for instance) where it 
may be applied in a
qualitative way.  (This echoes of course my catastrophe 
theory philosophy.)
Papers written in this state of mind are not read by 
professional
mathematicians, who see no need for communication with any 
other disciplines
apart from physics.  And they are not intelligible to 
people of the other
speciality, who generally lack the necessary mathematical 
culture.  As a result
they remain practically unread.  The case may be defended 
of papers which have
to create their own readership; they are babies without 
parents.
\vfill
\eject


\newletter

\letaddress
Edward Witten 
School of Natural Sciences
Institute for Advanced Study
Princeton, NJ 08540
IN\%"WITTEN@sns.ias.edu"
\endletaddress

Jaffe and Quinn attempt to comment on the role in 
mathematics of
some contemporary developments in physics.  I feel that 
the article
(in the section ``New Relations with Physics'')
conveys a rather limited idea of the role in physics of 
some of the
new developments in question.

Let me first try in one paragraph to summarize the state 
of knowledge
of physics.  (For a more extensive account, see the 
beginning of my
article on ``Physics and Geometry'' in the proceedings of 
the 1986
International Congress of Mathematicians.)
Gravitation is described at the classical level
by general relativity, which is based on Riemannian 
geometry.
Straightforward attempts at extending general relativity 
to a quantum
theory have always led to extremely severe difficulties.
Other observed forces are described by a quantum gauge 
theory (the ``standard
model''), whose construction  involves
(in addition to the machinery of quantum field theory)
the choice of a Yang-Mills gauge group ($SU(3)\times SU(2)
\times U(1)$ encompasses the known interactions); a 
representation of
that group for charged fermions (experiment indicates a 
rather complicated
reducible representation, related to phenomena such as 
parity violation
and the fractional electric
charges of quarks); and a relatively little understood 
mechanism of symmetry
breaking.

The main unsolved problems are generally considered to be to
overcome the inconsistency between gravity and quantum 
mechanics;
to unify the various other forces with each other and with 
gravity;
and to understand symmetry breaking and the vanishing of 
the cosmological
constant.

In the early 1980s,
it became clear---through the work of M. B. Green, J. H. 
Schwarz,
and L. Brink, building on pioneering contributions of 
others from the
1970s---that string theory
offered a framework (in my view the only promising
framework known) for overcoming the inconsistency between
gravity and quantum mechanics.  Actually, that is a 
serious understatement.
It is not just that in string theory, unlike previous 
frameworks of
physical theory, quantum gravity is possible; rather, the 
existence of
gravity is an unavoidable prediction of string theory.  In 
the early
development of the theory,
literally dozens of papers were written in an unsuccessful 
effort
to eliminate the features that lead to the prediction of 
gravity.

By the early 1980s, it was fairly clear
(from overwhelming circumstantial
evidence, not a mathematical theorem) that string theory 
made sense
and predicted gravity, but it appeared extremely difficult 
to
apply string theory to nature.  The reason for this was 
that at the time,
it appeared impossible in the context of string theory for 
the weak
interactions to violate parity.  Then in 1984, this 
difficulty
was overcome as a result of a new theoretical insight, and 
as a bonus the
gauge group and fermion representation of the standard 
model suddenly
emerged rather naturally from the theory.

On a more theoretical side, supersymmetry (or bose-fermi 
symmetry;
supergeometry) is another general prediction of string 
theory.  World-sheet
supersymmetry was invented by P. Ramond in 1970 to\ 
incorporate fermions in
string theory; fermions exist in nature, so this was 
necessary to make string
theory more realistic.  Space-time supersymmetry was 
invented by J. Wess and B.
Zumino in 1974 based on an analogy with world-sheet 
supersymmetry.  Ever since
then, supersymmetry has fascinated physicists, especially 
in connection to the
little-understood symmetry-breaking mechanism of the 
standard model.
Supersymmetry is not an established experimental fact, 
though a possible
partial explanation for the measured values of the strong, 
weak, and
electromagnetic coupling constants based on supersymmetry 
has attracted much
interest.  There is an active search for more direct 
experimental confirmation
of supersymmetry at high-energy accelerators; this is 
regarded by many as one
of the prime missions of the proposed Superconducting 
Supercollider.

The main immediate obstacle to progress in extracting
more detailed experimental predictions from string theory 
(beyond generalities
such as the existence of gravity) would appear to be that 
the vanishing
of the cosmological constant is not understood 
theoretically.

More fundamentally,
I believe that the main obstacle is that the core
geometrical ideas---which must underlie string theory the 
way Riemannian
geometry underlies general relativity---have not yet been 
unearthed.
At best we have been able to scratch the surface and 
uncover things
that will most probably eventually be seen as spinoffs of 
the
more central ideas.  The search for these more central ideas
is a ``mathematical'' problem which at present preoccupies 
primarily
physicists.  Some of the spinoffs have, however,
attracted mathematical interest in different areas.

In general, I think that the motivations for string theory
in physics are much stronger and more focussed than
Jaffe and Quinn convey.


\newletter

\letaddress
Sir Christopher Zeeman, FRS
Gresham Professor of Geometry
Hertford College
Oxford 0X1 3BW
England
\endletaddress

Dear Dick,

        Thank you for your invitation to respond to the 
paper by Jaffe \& Quinn
on Theoretical Mathematics [1].

        Their account of catastrophe theory is misleading, 
because Ren\'e
Thom's work on singularities [4] was firm.  In his 
subsequent development of
catastrophe theory he focused attention upon the key 
unsolved steps in the
underlying mathematics by making specific conjectures and 
encouraging
Malgrange, Mather and others to prove them.  For example 
Malgrange writes in
the introduction to his 1966 book [2] on differentiable 
functions:
\medskip
        {\parindent=20pt \narrower\noindent In particular, 
I consider it my
duty to state that one of the main results, ``the 
preparation theorem for
differentiable functions'', was proposed to me as a 
conjecture by R.~Thom, and
that he had to make a great effort to overcome my initial 
scepticism.\par}

        In my own book [7] I gave a complete and 
mathematically rigorous proof
of the classification of elementary catastrophes of 
codimension $\le 5$, and
the $C^\infty$-density of generic global parametrised 
functions.  I also made a
number of scientific models and scientific predictions, 
several of which have
been subsequently confirmed by experimentalists.  This 
does not fit the
description given by Jaffe and Quinn of being 
``mathematically theoretical''
(in their terminology).  In fact there have been hundreds 
of successful
scientific applications of catastrophe theory.

        What controversy there was about catastrophe 
theory was short-lived for
two reasons:  firstly the underlying mathematics was 
rigorous, and secondly the
critics were not scientists but a few journalists and 
mathematicians who were 
ignorant of the science and did not fully understand the 
mathematics.  
For example the scientific mistakes in [5] were answered 
in [6], and the 
mathematical mistakes in [3] were explained in [8].

        Turning to Jaffe and Quinn's main thesis, I 
applaud their appeal to
authors to distinguish more clearly between theorems and 
conjectures, and I
deplore their suggestion that the mathematical community 
should mimic the
physics community by separating those who make conjectures 
from those who prove
theorems.  The best mathematicians have always done both, 
and always will.

\lethead{References}

\item{1.} A. Jaffe and F. Quinn, {\it Theoretical 
Mathematics\RM{:} 
Toward a cultural 
        synthesis of mathematics and theoretial 
physics\/}, Bull. Amer. Math.
	Soc. {\bf 29} (1993), 1-13.

\item{2.} B. Malgrange, 
{\it Ideals of differentiable functions\/}, Oxford Univ. 
Press, 
Oxford, 1966.

\item{3.} S. Smale, {\it Review of E.C. Zeeman\RM{:}
 Catastrophe theory, selected papers\/}
        1972--1977, Bull. Amer. Math. Soc. {\bf 84} 
(1978), 1360-1368.

\item{4.}  R. Thom, {\it Les singularit\'es des applications
diff\'erentiables\/},
	Ann. Inst. Fourier {\bf 6} (1956), 43-87.

\item{5.}  R.S. Zahler and H. Sussmann, 
{\it Claims and accomplishments of applied
        catastrophe theory\/}, Nature {\bf 269} (1977), 
759-763.

\item{6.} {\it Correspondence on catastrophe theory\/},  
Nature {\bf 270} (1977), 
381--384
        and 658.

\item{7.}  E. C. Zeeman,  {\it Catastrophe theory, 
selected papers\/} 
1972--1977, 
        Addison Wesley, Reading, MA, 1977.

\item{8.} E. C. Zeeman, {\it Controversy in science\RM{:}
 On the ideas of Daniel Bernoulli
        and Ren\'e Thom\/},  
The 1992/3 Johann Bernoulli Lecture, Gr\"oningen (to
        appear in Nieuw Archief van de Wiskunde).

\enddocument